\documentclass[11pt]{amsart}
\usepackage[marginratio=1:1,height=600pt,width=400pt,tmargin=117pt]{geometry}
\usepackage{amsmath,amsthm,amssymb}
\usepackage{comment}
\usepackage{enumitem}
\usepackage{cleveref}
\usepackage{mathtools}
\usepackage{float}
\usepackage{tikz}
\usetikzlibrary{calc}
\usepackage[final]{pdfpages}

\newtheorem{theorem}{Theorem}[section]
\newtheorem{lemma}[theorem]{Lemma}
\newtheorem{proposition}[theorem]{Proposition}

\theoremstyle{definition}
\newtheorem{definition}[theorem]{Definition}

\theoremstyle{remark}
\newtheorem{remark}[theorem]{Remark}

\numberwithin{equation}{section}

\newcommand{\R}{\mathbb{R}}
\newcommand{\N}{\mathbb{N}}
\newcommand{\Z}{\mathbb{Z}}

\DeclareMathOperator{\diam}{diam}
\DeclareMathOperator{\dist}{dist}

\DeclareMathOperator{\Out}{Out}
\DeclareMathOperator{\inte}{int}

\makeatletter
\newsavebox\myboxA
\newsavebox\myboxB
\newlength\mylenA

\newcommand*\xoverline[2][0.75]{%
    \sbox{\myboxA}{$\m@th#2$}%
    \setbox\myboxB\null% Phantom box
    \ht\myboxB=\ht\myboxA%
    \dp\myboxB=\dp\myboxA%
    \wd\myboxB=#1\wd\myboxA% Scale phantom
    \sbox\myboxB{$\m@th\overline{\copy\myboxB}$}%  Overlined phantom
    \setlength\mylenA{\the\wd\myboxA}%   calc width diff
    \addtolength\mylenA{-\the\wd\myboxB}%
    \ifdim\wd\myboxB<\wd\myboxA%
       \rlap{\hskip 0.5\mylenA\usebox\myboxB}{\usebox\myboxA}%
    \else
        \hskip -0.5\mylenA\rlap{\usebox\myboxA}{\hskip 0.5\mylenA\usebox\myboxB}%
    \fi}
\makeatother

\title[The 2D Euler equations]
      {On stable solutions to the Euler equations in convex planar domains}

\author[Bian Wu]{Bian Wu}
\email{bian.wu@math.ethz.ch, bwu@ias.edu}
\address{Institute for Advanced Study , 1 Einstein Drive, Princeton, NJ 08540, USA}
\date{\today}

\begin{document}

\begin{abstract}
In convex planar domains, given an initial vorticity with one sign, we study the regularity and geometric properties of the dynamically stable solutions to the Euler equations in the coadjoint orbit of the initial vorticity. These flows have elliptic stagnation points. Under some nondegeneracy conditions on the data, we show they are H\"older continuous and have convex level curves. We also give a detailed description for the set of stagnation points. If the initial vorticity has nice level set topology, these stable solutions are in the $L^\infty$-strong closure of the coadjoint orbit. We also demonstrate the sharpness of most assumptions we made.
\end{abstract}

\maketitle

\section{Introduction}

The 2D Euler equations in a bounded domain $\Omega \subset \R^2, \, \partial \Omega \in C^1$ with initial condition $\omega(\cdot,0) = \omega_0 \in L^{\infty}(\Omega)$ are given by
\begin{equation} \label{evo_euler_vor}
  \begin{split}
  \partial_t \omega + u \cdot \nabla \omega &= 0, \\
  u &= \nabla^\perp \Delta^{-1} \omega,\\
  u \cdot n &= 0, \quad\quad \text{at } \partial \Omega.
  \end{split}
\end{equation}
Here $n$ is the outward normal vector on $\partial \Omega$. The vorticity $\omega = \partial_1 u_2 - \partial_2 u_1$ is conserved along particle trajectory. Indeed, assuming enough regularity on $\omega_0$, if we denote by $\Phi:\Omega \times [0,\infty) \rightarrow \Omega$ the flow map corresponding to the velocity field $u$, then $\Phi$ satisfies
\[ \partial_t \Phi(x,t) = u( \Phi(x,t),t ), \quad \Phi(x,0)=x \]
and $\omega(x,0)=\omega(\Phi(x,t),t)$ for any $t>0$. If $\omega_0$ is smooth, $\{\Phi(\cdot,t)\}_{t \geq 0}$ is a family of volume-preserving diffeomorphisms. Define the coadjoint orbit $\mathcal{O}_{\omega_0}$ of a vorticity $\omega_0$ to be
\begin{equation} \label{orbit_def}
  \mathcal{O}_{\omega_0} := \{ \omega_0 \circ \eta^{-1} : \eta \in \mathcal{D}_{\text{vol}}(\Omega)\},
\end{equation}
where $\mathcal{D}_{\text{vol}}(\Omega)$ is the group of volume-preserving diffeomorphisms from $\Omega$ to $\Omega$. Then 2D Euler equations with initial vorticity $\omega_0$ can be considered as a dynamical system in $\mathcal{O}_{\omega_0}$.
\par

Since the velocity field $u$ is divergence-free and $u \cdot n = 0$, there exists a stream function $\psi : \Omega \rightarrow \R$ with $u = (-\partial_2 \psi, \partial_1 \psi)^T := \nabla^\perp \psi$ and $\omega = \Delta \psi$. The stationary 2D Euler equations can be written as
\begin{equation} \label{sta_euler}
  \begin{split}
  \{\psi, \omega\} &= 0, \\
  \omega &= \Delta \psi,
  \end{split}
\end{equation}
where $\{\cdot,\cdot\}$ is the 2D Poisson bracket defined by $\{f,g\} := \partial_1f \partial_2g - \partial_2f \partial_1g$. Kelvin \cite{thomson1910} pointed out that the stationary 2D Euler equations \eqref{sta_euler} have a nice variational structure. Namely, if $\psi$ is a critical point of the kinetic energy functional $E$ in an orbit $\mathcal{O}_{\omega_0}$, where $E$ is defined by
\begin{align} \label{kinetic_energy}
  E(\omega) := \frac{1}{2} \int_\Omega |\nabla \psi|^2 dx,
\end{align}
$\psi$ solves the stationary 2D Euler equations \eqref{sta_euler}. This functional is formally quadratic and corresponds to the kinetic energy of the velocity field. Therefore, it is natural to ask of the existence and the regularity of minimizers, see Remark 2.9 in Chapter 2 in Arnold and Khesin. The existence is known to fail in $C^\infty$ topology in general because of the constraints on level set topology showed by Ginzburg and Khesin in \cite{MR1288220}, or because of regularity issues mentioned by Nadirashvili \cite{MR3067825} and also proved in \Cref{counterexample_naive}. On the other hand, the existence holds in weaker topologies, such as weak-$*$ topology in $L^\infty$ showed by Shnirelman \cite{MR1240495} or strong topology in $L^q, \, q \in [1,\infty)$ by Burton \cite{MR870963} and \cite{MR998605}. Burton also showed that $E$ has infinitely many critical points in the strong closure of $\mathcal{O}_{\omega_0}$ in $L^q$ and the minimizer is dynamically stable in \cite{MR1135975} and \cite{MR2186035}.

\begin{theorem}[Burton, \cite{MR1135975}, \cite{MR2186035}] \label{burton_thm}
Given a bounded open domain $\Omega \subset \R^2$ with $1 \leq q < \infty$. Let $\omega_0 \in L^q(\Omega)$ satisfies $\omega_0 \geq 0$ or $\omega_0 \leq 0$, then the functional $E$ defined in \eqref{kinetic_energy} admits a unique minimizer, at least one maximizer, and infinitely many saddle points of $E$ in the $L^q$-strong closure of $\mathcal{O}_{\omega_0}$. Moreover, the minimizer $\omega$ is dynamically stable in the $L^q$-topology of the vorticity and it satisfies $\omega = f(\Delta^{-1} \omega)$ for some nondecreasing function $f$.
\end{theorem}

In this paper we establish the existence of minimizer in the strong $L^\infty$-closure of $\mathcal{O}_{\omega_0}$ under some nondegeneracy conditions. Note that $L^\infty$ is the natural space for the wellposedness of 2D Euler equations. We show the global minimizer is H\"older continuous under these conditions.

\begin{theorem} \label{main_thm_intro}
Given a bounded convex domain $\Omega \subset \R^2$ with $C^{2,\alpha/4}$ boundary, $\alpha \in (0,1]$. Let $\omega_0 \in C^\alpha(\bar{\Omega})$ and $\omega_0 \geq 0$.
\begin{enumerate}[leftmargin=*,label=\textup{(\arabic*)},align=left]
    \item If $\inf_{x \in \Omega} \omega_0(x) >0$, then the global minimizer $\omega$ of the functional $E$ in \Cref{burton_thm} is in $C^{\alpha/4}(\bar{\Omega})$. Also, the set of stagnation points is $\{\psi = \min \psi\}$, and it contains a single point or a line segment.
    \item Assume $\omega_0$ is constant on $\partial \Omega$ and the set $\{\omega_0 < s\}$ is simply connected for any $s \in (\inf \omega_0, \sup \omega_0)$, then $E$ has a unique global minimizer in the $L^\infty$-strong closure of $\mathcal{O}_{\omega_0}$ which coincides with $\omega$ in \Cref{burton_thm}.
    \item Assume $\omega_0$ is not constant on $\partial \Omega$ or the set $\{\omega_0 < s\}$ is not simply connected for some $s \in (\inf \omega_0, \sup \omega_0)$, then $E$ has no minimizer in the $L^\infty$-strong closure of $\mathcal{O}_{\omega_0}$.
    \item Suppose $\inf_{x \in \Omega} \omega_0(x) >0$. Assume $\omega_0$ is constant on $\partial \Omega$ and has a single critical point, then $\omega$ is smooth in any compact subset of $\Omega \backslash \{\psi = \min \psi\}$.
\end{enumerate}
\end{theorem}

Secondly, we show that dynamically stable solutions in convex domains have convex level curves.

\begin{theorem} \label{thm_convexity}
In a bounded convex domain $\Omega \subset \R^2$, suppose that $\psi \in W^{2,2}(\Omega)$ solves the steady Euler equations \eqref{sta_euler} and satisfies Arnold's stability criterion weakly. Then $\psi$ has convex level curves. Here, Arnold's stability criterion are defined in \Cref{arnold_sta} and \Cref{arnold_sta_weak}.
\end{theorem}

\begin{remark}
In \Cref{main_thm_intro}, we have all conclusions symmetrically for $\sup_{x \in \Omega} \omega_0(x) <0$ or $\omega_0 \leq 0$. Also, (2) and (3) in \Cref{main_thm_intro} hold for the maximizers which satisfy Arnold's stability criterion weakly. If we further assume the vorticity of the maximizer is not zero in $\bar{\Omega}$, (1) and (4) in \Cref{main_thm_intro} hold. Therefore, \textit{we have all nice properties in \Cref{main_thm_intro} and \Cref{thm_convexity} for dynamically stable solutions.} On the contrary, it is easy to see the H\"older regularity fails in general for the saddle points in \Cref{burton_thm}. Nadirashvili \cite{MR3483474} also proved the convexity of the level curves fails in general for the saddle points.
\end{remark}

\begin{remark}
(2) and (3) shows the constraints on the level set topology is necessary for the existence of the minimizer in the $L^\infty$-strong closure of $\mathcal{O}_{\omega_0}$. The continuity of $\omega_0$ is also necessary. A counterexample about vortex patch is given in \Cref{counterexample_3}.
\end{remark}

\begin{remark} \label{rmk_cont}
Although we assume $\omega_0 \in C^\alpha(\bar{\Omega})$, general vorticities in the $L^\infty$-strong closure of $\mathcal{O}_{\omega_0}$ may not belong to $C^\beta$ for any $\beta>0$. In the $L^\infty$ weak-$*$ closure or $L^q$-strong closure of $\mathcal{O}_{\omega_0}$, $q \in [1,\infty)$, the vorticities are only bounded. \Cref{main_thm_intro} says that the stable solutions have better regularity than general vorticities in these closures of $\mathcal{O}_{\omega_0}$.
\end{remark}

\begin{remark}
We impose nondegeneracy conditions directly depending on $\Omega$, namely the assumption about the convexity. The assumption on the sign of $\omega_0$ is necessary for the existence of the global minimizer. Indeed, Burton and McLeod \cite{MR1135975} proved that the global minimizer does not exist when $\omega_0$ admits both positive values and negative values. These conditions ensure that we are only dealing with elliptic stagnation points. When either of these fails, there might be hyperbolic stagnation points which genuinely complicates the picture. In (1) and (4), $\inf_{x \in \Omega} \omega_0(x) = 0$ causes some degeneracy which our argument cannot handle. The assumption on the regularity of $\partial \Omega$ is needed even for the Schauder theory for linear elliptic equations.
\end{remark}

\begin{remark}
Under the assumptions in (4) of \Cref{main_thm_intro}, the smoothness does not hold at $\{\psi = \min \psi\}$ in general. We believe this was already known, at least to Nadirashvili \cite{MR3067825}. Since Nadirashvili did not give the details, we give a concise counterexample to the smoothness in \Cref{counterexample_naive} of \Cref{counterexample_naive_sec} for the sake of completeness.
\end{remark}

In \Cref{main_thm_intro}, we study the elliptic stagnation points of $\psi$. This case was left over by several previous works. Nadirashvili \cite{MR3067825} showed that, in smoothly bounded convex domain $\Omega$ with $\omega_0 > 0$ having a single critical point, away from the stagnation point, the global minimizer is smooth and admits analytic level curves. In the case of annulus, assuming $\psi = \Delta^{-1} \omega$ has no stagnation point and some other nondegeneracy conditions, Choffrut and Sverak \cite{MR2899685} showed, nearby a steady-state, there is a local one-to-one correspondence between smooth steady states and coadjoint orbits. Therefore, the flows with elliptic stagnation points are unclear. In \Cref{main_thm_intro}, we prove that $u = \nabla^\perp \psi \in C^{1,\alpha}(\bar{\Omega})$, and Arnold stable solutions have convex level curves. This confirms Nadirashvili's guess on page 732 and page 733 of \cite{MR3067825}.
\par

The coherent vortex structures in 2D turbulence suggests that the final states of 2D Euler equations are more regular than $L^\infty$ in large scales, see \cite{MR2929420}, \cite{MR1682816} and \cite{MR1121395}. This has been supported by numerical and physical studies, for example, \cite{clercx1998spontaneous}, \cite{matthaeus1991decaying}, \cite{MR2449898} and \cite{MR1962736}. We also refer to \cite{dolce2022maximally} and \cite{drivas2022singularity} for more mathematical discussions, to \cite{MR3355006} and \cite{montgomery1992relaxation} for the relations with inviscid limits. As mentioned in \Cref{rmk_cont}, \Cref{main_thm_intro} tells that the dynamically stable solutions (also termed as weak attractors in some literature) have more regularity then general vorticities in suitable closure of $\mathcal{O}_{\omega_0}$. This is consistent with relevant physical predictions.
\par

The key tool in this work is the following. The stream function for the minimizer $\psi:\Omega \rightarrow \R$ satisfies
\begin{equation} \label{semilinear_Euler}
  \begin{split}
    \Delta \psi &= f(\psi), \\
    \psi|_{\partial \Omega} &= 0
  \end{split}
\end{equation}
for some nondecreasing function $f:\R \rightarrow \R$. The difference of this work from other works on semilinear elliptic equations is that we only have very limited quantitative information on $f$. Most techniques developed for smooth $f$ with nice properties cannot be applied directly. Our main task is to deduce the quantitative information on $f$ and on $\psi$ based on the formal indentity $f = (\omega^*)^{-1} \circ \psi^*$ and the monotonicity of $f$. \par

The rest of the paper is organized as follows. In \Cref{sec_stability}, we give some preliminaries on Arnold's stability criterion. In \Cref{level_set}, we prove that the level sets of $\psi$ are convex curves. In \Cref{convex_geometric_con}, we prove a geometric lemma for convex rings. In \Cref{regularity}, we use these two ingredients to show the H\"older continuity of $\psi^*$ and $(\omega^*)^{-1}$. Here, $\psi^*$ is defined by
\begin{align} \label{dis_eq}
  \psi^*(t) := |\{x \in \Omega \mid \psi(x) \leq t \}|
\end{align}
and $(\omega^*)^{-1}$ is a left-inverse of $\omega^*$ defined in \Cref{diff_lemma_1}. Then we use $f = (\omega^*)^{-1} \circ \psi^*$ to prove H\"older continuity in \Cref{main_thm_intro}. In \Cref{critical_points}, we prove $\{\psi = \min \psi\}$ is exactly the set of stagnation points, namely (1) and (4) in \Cref{main_thm_intro}. In \Cref{counterexamples}, we prove (2) and (3) in \Cref{main_thm_intro} and give a counterexample to show the sharpness of the conditions. In \Cref{counterexample_naive_sec}, we give a counterexample to show the smoothness may fail on the set $\{\psi = \min \psi\}$.

\vspace{3mm}

\noindent \textbf{Notation.} We use $\Delta^{-1}$ to denote the operator solving Laplace equation with Dirichlet boundary condition. A \textit{convex curve} is a set which is the boundary of some nonempty bounded closed convex set in $\R^2$. We use $|\cdot|$ to denote the Lebesgue measure in $\R^2$. We use $\diam A$ to denote the diameter of a bounded set $A \subset \R^2$, namely $\diam A := \sup_{x,y \in A} |x-y|$. For $x,y \in \R^2$ with $x \neq y$, we use $\overline{xy}$ to denote the line segment connecting $x$ and $y$, i.e. $\overline{xy} := \{tx+(1-t)y \mid t \in [0,1]\}$. For $a \in \R$, we use $\lceil a \rceil$ to denote the smallest integer in the set $\{b \in \R \mid b \geq a\}$. We define the distance function $\dist (A,B) := \inf_{x \in A, y \in B} |x-y|$ between two sets $A,B \subset \R^2$. For a function $\psi: A \rightarrow \R$ with $A \subset \R^2$ and $|A| < \infty$, we define its distribution function to be \eqref{dis_eq}.

% !TEX root = P5_paper1.tex

\section{Preliminaries on Arnold stability criterion} \label{sec_stability}

In this section, we give definitions for Arnold's stability criterion and its weak form. We compare these two criterion in \Cref{sta_equivalence}. For more discussions, we refer to \cite{MR4268535} and \cite{drivas2022singularity} for Arnold's stability criterion, or to \cite{MR1245492} and \cite{MR2186035} for weaker alternatives of Arnold's stability criterion.

\begin{definition} \label{arnold_sta}
In a bounded domain $\Omega$ with $\partial \Omega \in C^{2,0+}$, a steady flow satisfies Arnold's stability criterion if the stream function $\psi \in W^{3,1}(\Omega)$ of this steady flow satisfies 
\begin{align} \label{arnold_eq0}
    C^{-1}|\nabla \psi| \leq |\nabla \Delta \psi| \leq C|\nabla \psi|
\end{align}
for some $C>0$ and either of the following conditions:
\begin{enumerate}[leftmargin=*,label=\textup{(\arabic*)},align=left]
    \item $\nabla \psi \cdot \nabla \Delta \psi = |\nabla \psi| |\nabla \Delta \psi|$.
    \item $-\nabla \psi \cdot \nabla \Delta \psi = |\nabla \psi| |\nabla \Delta \psi|$. And for the nonincreasing and Lipschitz function $f$ with $\omega = f(\psi)$, we have $\inf_{s \in [\inf \psi, \sup \psi]} f'(s) > -\lambda_1(-\Delta,\Omega)$. Here, $\lambda_1(-\Delta,\Omega)$ is the first eigenvalue of the Laplace operator.
\end{enumerate}
\end{definition}

\begin{definition} \label{arnold_sta_weak}
For a bounded domain $\Omega$, a steady flow in $\Omega$ satisfies Arnold's stability criterion \textit{weakly} if the stream function $\psi \in W^{2,2}(\Omega)$ of this steady flow satisfies either of the following conditions:
\begin{enumerate}[leftmargin=*,label=\textup{(\arabic*)},align=left]
    \item There exists a nondecreasing function $f \in L^\infty(\R)$ with $f \geq 0$ or $f \leq 0$, such that $\Delta \psi = f(\psi)$.
    \item There exists a nonincreasing function $f \in W_{\text{loc}}^{1,1}(\R) \cap L^\infty(\R)$ with $f \geq 0$ or $f \leq 0$, such that $\Delta \psi = f(\psi)$ and $\inf_{s \in [\inf \psi, \sup \psi]} f'(s) > -\lambda_1(-\Delta,\Omega)$.
\end{enumerate}
\end{definition}

\begin{remark}
Now we justify the existence of $f$ in \Cref{arnold_sta}. The regularity assumption $\psi \in W^{3,1}(\Omega)$ is the minimal regularity to make sense of Arnold's stability criterion in \cite{MR4268535}. Under this assumption, one can show $\omega = \Delta \psi$ is Lipschitz continuous. Indeed, $\psi \in W^{3,1}(\Omega)$ implies $\nabla \psi \in L^q$ for any $q \in [1,\infty)$, so $\nabla \omega \in L^q$ by \eqref{arnold_eq0}. Then $\omega \in C^\alpha$ for some $\alpha \in (0,1)$ and $\nabla \psi \in C^{1,\alpha}$, hence $\omega$ is Lipschitz by \eqref{arnold_eq0}. With $\psi \in C^{2,\alpha}$ and $\omega$ Lipschitz, one can conclude there exists a continuous function $f:\R \rightarrow \R$ with $\omega = f(\psi)$. One can also show $f$ is monotone and Lipschitz from the functional relation between $\omega$ and $\psi$. In (1), $f$ is nondecreasing. In (2), $f$ is nonincreasing.
\end{remark}

\begin{remark}
The difference between \Cref{arnold_sta} and \Cref{arnold_sta_weak} is the regularity of $f$. For the minimizer given by \Cref{burton_thm}, we know only it satisfies Arnold's stability criterion weakly. In general, for a steady flow which satisfies Arnold's stability criterion weakly, $f$ is not Lipschitz continuous. In \Cref{arnold_sta_weak} and \Cref{levelset_conv} below, one can allow $f \in L_{\text{loc}}^q$ for $q>1$. The conclusions of \eqref{levelset_conv} and \Cref{thm_convexity} hold after suitable modifications on our argument.
\end{remark}

 A steady flow satisfies Arnold's stability criterion implies its vorticity cannot have two signs. We prove this fact in the following lemma. The argument is based on maximum principle of the operator $-\Delta + f'$. Similar arguments will be needed in the next section. We omit the details there.

\begin{lemma} \label{sta_equivalence}
If a steady flow $\omega$ satisfies Arnold's stability criterion, then $\omega \geq 0$ or $\omega \leq 0$.
\end{lemma}

\begin{proof}
Assume $\omega$ admits two signs. Consider the case where (2) and \eqref{arnold_eq0} in \Cref{arnold_sta} are satisfied. Then there exists $c \in \R$ such that $f(c) = 0$ and $\psi = c$ is not empty. Without loss of generality, assume $c>0$. Consider the function $\varphi := \min\{c, \psi\} - c$. We know $\varphi \in H_0^1(\Omega)$ and $\varphi \neq 0$. Multiply $\Delta \psi = f(\psi)$ with $\varphi$ and integrate in $\Omega$, then we get
\begin{align*}
    \int_\Omega |\nabla \varphi|^2 dx = -\int_\Omega \varphi f(\varphi+c) dx.
\end{align*}
Since $f(\varphi+c) \in (-\gamma \varphi,0]$ for some constant $\gamma \in (0,\lambda_1(-\Delta,\Omega))$. Then we have
\begin{align*}
    \lambda_1(-\Delta,\Omega) > \gamma \geq \frac{\int_\Omega |\nabla \varphi|^2 dx} {\int_\Omega |\varphi|^2 dx}.
\end{align*}
This gives a contradiction, since $\lambda_1(-\Delta,\Omega)$ is the first eigenvalue of the Laplace operator in $\Omega$. The case where (1) and \eqref{arnold_eq0} in \Cref{arnold_sta} are satisfied follows from a similar argument or a direct application of maximum principle.
\end{proof}

\section{Level sets of the stream function} \label{level_set}

In this section we study the level curves of the stream function of weakly Arnold stable solutions.

\begin{proposition}   \label{levelset_conv}
For bounded convex planar domain $\Omega$, let $f:\R \rightarrow \R$ be bounded and non-decreasing function with $f \geq 0$. Suppose $f(s_0) = 0$ for some $s_0 \in \R$ Then the level sets $\{\psi = c\}$ of the solution $\psi$ to \eqref{semilinear_Euler} are convex curves for any $c \in (\min \psi, \max \psi]$, and $\{\psi = \min \psi\}$ is a convex set. Moreover, if $\inf_{s \in [\inf \psi, \sup \psi]} f(s)>0$, $\{\psi = \min \psi\}$ is a single point or a line segment.
\end{proposition}

\begin{remark}
If $\psi$ satisfies Arnold's stability criterion, namely $f$ is Lipschitz, $\{\psi = \min \psi\}$ is a single point. This is proved in \cite{MR3067825}.
\end{remark}

\Cref{thm_convexity} follows from \Cref{levelset_conv}. For the solutions where $f \leq 0$ or (2) in \Cref{arnold_sta_weak} is satisfied, the proof is analogues to \Cref{levelset_conv}. The convexity of level curves is a result of the following theorem on the convexity properties for elliptic equations. There is a vast literature on this topic. Here, we mention \cite{MR678504}, \cite{MR792181} and \cite{MR1623694} as references and we cite the following theorem.

\begin{theorem}[Cabr\'{e} and Chanillo \cite{MR1623694}, Finn \cite{MR2425011}]   \label{convexity_cite}
For bounded and strictly convex planar domain $\Omega$, let $f$ be a bounded smooth function satisfying $f \geq 0$. Suppose $f' \geq 0$ or $f' \in (-\lambda_1(-\Delta,\Omega), 0]$. Then the level sets of the solution of \eqref{semilinear_Euler} are convex curves. Here, being strictly convex means $\partial \Omega$ does not contain line segment.
\end{theorem}

Next we prove the convexity of the level curves for $f$ which is only bounded, non-decreasing and non-negative and $\Omega$ is not necessarily strongly convex.

\begin{proof}[Proof of \Cref{levelset_conv}]
We prove this lemma by approximation arguments. We need two levels of approximation, for $f$ and for $\Omega$.
\par

On the first level, we do the approximation for $f$. Fix $0 \leq \chi \in C_c^\infty([0,1])$ with $\int_{\R} \chi(s) ds=1$. Define a smoothing kernel $\{\chi_k\}_{k \geq 1}$ with $\chi_k(s) := k\chi(ks)$, then $f*\chi_k$ converges to $f$ locally in $L^q$ for any $q \in [1,\infty)$. Consider the function sequence $\{f_k\}_{k \geq 1}$ defined by $f_k(t) := (f * \chi_k)(t)$. Therefore, we have
\begin{enumerate}[leftmargin=*,label=\textup{(\arabic*)},align=left]
    \item $f_k \geq f$ for any $k \geq 1$.
    \item $f_k$ converges to $f$ locally in $L^q$ for any $q \in [1,\infty)$.
    \item $f_k(s) = 0$ for any $s \leq s_0-1$.
\end{enumerate}
Classical methods give the existence of a solution $\psi_k$ to \eqref{semilinear_Euler} with $f_k$. Here, we sketch a proof by calculus of variations. Define $F_k(t) = \int_{s_0-1}^t f_k(s)ds$ and consider the following functional
\[
    \frac{1}{2} \int_\Omega |\nabla \varphi|^2 dx + \int_\Omega F_k(\varphi) dx, \quad \varphi \in H_0^1(\Omega).
\]
We have that $F_k$ is positive, increasing and Lipschitz, thus it is easy to show this functional is coercive and weakly sequentially lower semi-continuous in $H_0^1(\Omega)$. Then it has a minimizer $\psi_k$ and it satisfies $\Delta \psi_k = f_k(\psi_k)$ weakly. Moreover, since $f_k$ is uniformly bounded in $L^\infty$, $\{\psi_k\}_{k \geq 1}$ is uniformly bounded in $C^{1,\beta}$ for some $\beta>0$ and it converges to a function $\psi_0$ in $C^{1,\beta-}$.
By a maximum principle argument, for non-decreasing $f$, the solution to \eqref{semilinear_Euler} is unique, so $\psi=\psi_0$. Similarly by a maximum principle argument, we also have $\psi_k \leq \psi$ in $\Omega$. \par

Since $f_k$ satisfies all conditions in \Cref{convexity_cite}, then the sublevel sets $\{\psi_k \leq c\}$ of $\psi_k$ are convex for all $k$ and all $c \in (\min \psi_k, \max \psi_k)$. Consider the set
\begin{align*}
    S_c := \{ x \in \Omega \mid \exists \{x_k\}_{k \geq 1} \subset \Omega
        \text{ with } \psi_k(x_k) \leq c \text{ and } \lim_{k \rightarrow \infty} x_k = x \}.
\end{align*}
By the $C^{1,\beta-}$ convergence of $\{\psi_k\}_{k \geq 1}$, we deduce that $S_c$ is a convex set and $S_c \subset \{\psi \leq c\}$. Since $\psi > \psi_k$ for any $k \geq 1$, $\{\psi \leq c\} \subset S_c$. Then the sublevel sets $\{\psi \leq c\}$ of $\psi$ are convex. \par

By strong maximum principle, $\{\psi = 0\}$ is a convex curve. Fix any $c \in (\min \psi, \max \psi)$ and consider the set $\bigcup_{t < c} \{\psi \leq t\}$. If $\{\psi \leq c\} \backslash \overline{\bigcup_{t < c} \{\psi \leq t\}}$ is empty, then the level set $\{\psi = c\}$ is a convex curve. If $\{\psi \leq c\} \backslash \overline{\bigcup_{t < c} \{\psi \leq t\}}$ is not empty, then $\{\psi=c\}$ has positive measure. In the interior of $\{\psi=c\}$, we have $\Delta \psi = 0$. By the assumptions on $f$, we can deduce $c=\min \psi$ and $\{\psi \leq \min \psi\}$ is convex. This proves the desired statements for smoothly bounded and strictly convex planar domain $\Omega$. \par

Now we sketch the approximation procedure for not necessarily strongly convex $\Omega$. Namely, we approximate $\Omega$ by bounded and strongly convex $\Omega_k$ with $\Omega \subset \Omega_k$ for any $k \geq 1$. This is similar to designing $f_k$. Then we look for the solution to $\Delta \psi_k = f_k(\psi_k)$ in $\Omega_k$ and repeat the argument above. This concludes our proof.

\end{proof}

\section{A lemma in convex geometry} \label{convex_geometric_con}

In this section, we prove a geometric lemma to answer the following question: Given nonempty bounded convex sets $A,D \subset \Omega$ with $\bar{D} \subset A$, what is the relation between the radius of the biggest ball enclosed in the convex ring $A\backslash \bar{D}$ and the area of the convex ring $A\backslash \bar{D}$? The answer is that we have a lower bound for the radius given in the following proposition.

\begin{proposition} \label{geometry_lem_ann}
There exists an absolute constant $\epsilon>0$ satisfying the following condition: Given any open convex sets $A, D \subset \R^2$ with and $\bar{D} \subset A$, there exists $B_r(x) \subset A\backslash \bar{D}$ for some $x \in A\backslash \bar{D}$ and some $$r \geq \frac{\epsilon|A\backslash \bar{D}|}{\diam D}.$$
\end{proposition}

To prove \Cref{geometry_lem_ann}, we need a generalized notion for outer nomal vector on Lipschitz graphs and a lemma.

\begin{definition}[Generalized outer nomal vector]
For any open set $D$ such that $\partial D$ is a Jordan curve given by a injective continuous map from $S^1$ to $\R^2$ and any $x \in \partial D$, we define the outer normal direction of $x$ as follows
\begin{align*}
    \Out_{\partial D}(x) := \{ &\theta \in \R / 2\pi \Z \mid \exists B_r(y) \subset \R^2 \backslash D\text{ for some } y \in \R^2, r > 0, \\
    &\text{such that } x \in \partial B_r(y) \text{ and } y-x=(r\cos \theta, r\sin \theta) \}.
\end{align*}
\end{definition}

\begin{lemma} \label{basic_convex}
For any bounded convex set $A \subset \R^2$ and $r>0$, we define
\begin{align*}
    N_r(A) := \{ x \in \R^2 \,|\, x \notin A, \dist(x,A) \leq r \}.
\end{align*}
Then 
$$ |N_r(A)| \leq 2\pi r\diam A + \pi r^2. $$
\end{lemma}

\Cref{basic_convex} is an intuitive analogue to the construction of the Riemann integral in elementary calculus, modulo some standard approximation argument. We omit its proof here. Now we are in the position to prove \Cref{geometry_lem_ann}. The idea is as follows: Let $r>0$ be the radius of the biggest ball enclosed in $A \backslash \bar{D}$, then $N_r(D)$ covers most area of $A \backslash \bar{D}$. The exceptional area $(A \backslash \bar{D}) \backslash N_r(D)$ is small and can be estimated using the convexity.

\begin{proof}[Proof of \Cref{geometry_lem_ann}]
Define $$ R:= \sup\{ r>0 \mid B_r(x) \subset A\backslash\bar{D} \text{ for some } x \in A\backslash\bar{D}\}. $$
Let $G := (A\backslash \bar{D}) \backslash N_R(D)$, then
$$ |A\backslash \bar{D}| \leq |N_R(D)| + |G|. $$
The main task of this proof is to estimate $|G|$. For this purpose, we introduce the following notions. For any point $x \in \partial D$ and any $\theta \in \Out_{\partial D}(x)$, define a function $\rho: \partial D \times \R / 2\pi \Z \rightarrow \R^+$ as
\begin{align*}
    \rho(x,\theta) = \sup \{ &r>0 \mid \exists B_r(y) \subset A \backslash \bar{D}\text{ for some } y \in \R^2, \\
    &\text{such that } x \in \partial B_r(y) \text{ and } y-x=(r\cos \theta, r\sin \theta) \}.
\end{align*}
We also define a function $z: \partial D \times \R / 2\pi \Z \rightarrow \R^2$ such that $B_{\rho(x,\theta)}(z(x,\theta)) \subset A \backslash \bar{D}$ and $z(x,\theta)-x=(r\cos \theta, r\sin \theta)$. Define a set $H$ as
$$ H:= \bigcup_{x \in \partial D, \theta \in \Out_{\partial D}(x)} B_{\rho(x,\theta)}(z(x,\theta)). $$
Denote $J:=(A\backslash \bar{D}) \backslash \bar{H}$. Note that $H \subset N_R(D)$, hence
\begin{align} \label{geometry_lem_ann_eq4}
    |A\backslash \bar{D}| \leq |N_R(D)| + |J|.
\end{align}
Because $D$ is a convex set and $\bar{D} \subset A$, we have $\partial J \cap \partial D = \varnothing$. From the definition of $H$ and the convexity of $A$ and $D$, we know each connected component of $J$ is simply connected and its boundary intersects $\partial A$.

We summarize the following facts (see Figure \ref{fig2} for details):
\begin{enumerate}[leftmargin=*,label=\textup{(\arabic*)},align=left]
    \item $J$ has at most countably many connected components denoted by $\{J_k\}_{k \in \Lambda}$. Here $\Lambda$ is a countable set or a finite set.
    \item For each connected component $J_k$, $\partial A \cap \partial J_k$ is homeomorphic to $[0,1] \subset \R$. To fix the orientation, we denote the two endpoints of $\partial A \cap \partial J_k$ by $x_k^1, x_k^2$ such that $\inf \Out_{\partial A}(x_k^2) - \sup \Out_{\partial A}g(x_k^1) \in (0,\pi)$ or $\inf \Out_{\partial A}(x_k^2) - \sup \Out_{\partial A}g(x_k^1) \in (-2\pi,-\pi)$.
    \item For each connected component $J_k$, there exists $x_k \in \partial A$ such that for some $\theta_k \in \Out_{\partial A}(x_k)$, such that $x_k^1, x_k^2 \in \partial B_{\rho(x_k,\theta_k)}(z(x_k,\theta_k))$.
    \item Define a sequence $\{\zeta_k\}_{k \in \Lambda}$ by
    \begin{align*}
         \zeta_k = 
            \begin{cases}
            \inf \Out_{\partial A}(x_k^2) - \sup \Out_{\partial A}g(x_k^1), &\text{if }\inf \Out_{\partial A}(x_k^2) - \sup \Out_{\partial A}g(x_k^1) \in (0,\pi),\\
            \inf \Out_{\partial A}(x_k^2) - \sup \Out_{\partial A}g(x_k^1) + 2\pi, &\text{if }\inf \Out_{\partial A}(x_k^2) - \sup \Out_{\partial A}g(x_k^1) \in (-2\pi,-\pi).
            \end{cases}
    \end{align*}
    We have $\sum_{k \in \Lambda} \zeta_k \leq 2\pi$.
    \item $R= \sup_{x \in \partial D, \theta \in \Out_{\partial D}(x) } \rho(x,\theta)$.
\end{enumerate}

We only prove the third fact by contradiction. Others are easy to verify. Assume there exists no $x_k \in \partial A$ with some $\theta_k \in \Out_{\partial A}(x_k)$, such that $x_k^1, x_k^2 \in \partial B_{\rho(x_k,\theta_k)}(z(x_k,\theta_k))$. Then as in Figure \ref{fig1}, by the definitions of $H$ and $J$, $x_k^1$ and $x_k^2$ must lie on two distinct balls with centers $y_1:=z(y_k^1,\theta_k^1)$ and $y_2:=z(y_k^2,\theta_k^2)$ and with radius $\rho(y_k^1,\theta_k^1)$ and $\rho(y_k^2,\theta_k^2)$ respectively, for some $y_k^1, y_k^2 \in \partial D$, $\theta_k^1 \in \Out_{\partial D}(y_k^1)$, $\theta_k^2 \in \Out_{\partial D}(y_k^2)$ with $y_k^1 \neq y_k^2$. Now we can deform from one ball to another ball continuously in $A \backslash \bar{D}$ while tangentially touching $\partial A$ and $\partial D$. Then there exists $x_k \in \partial A$ with some $\theta_k \in \Out_{\partial A}(x_k)$ such that $B_{\rho(x_k,\theta_k)}(z(x_k,\theta_k))$ intersects with $\partial J_k$ at another two points. Moreover, we can deduce $x_k^1, x_k^2 \notin \partial J_k$. This gives a contradiction.

\begin{figure}[H]
\centering
\begin{tikzpicture}
  [
    scale=0.5,
    %>=stealth,
    point/.style = {draw, circle,  fill = black, inner sep = 1pt},
    %dot/.style   = {draw, circle,  fill = black, inner sep = .2pt},
  ]

  % the circle
  \def\rad{2.4}
  \node (o1) at (0,0) [point, label = {below:$z(x_k,\theta_k)$}]{};
  \draw (o1) circle (\rad);

  \node (o2) at (3.6,-1.2) [point, label = {below:$y_2$}]{};
  \draw (o2) circle (1.2);

  \node (o3) at (-4.4,-1.8) [point, label = {below:$y_1$}]{};
  \draw (o3) circle (1.2);

  % triangle nodes: just points on the circle
  \node (n1) at (0,3) [label = {left:$\partial J_k$}]{};
  \node (n2) at (-8,-3) {};
  \node (n3) at (8,-3) [label = {below right:$\partial A$}] {};

  \node (n4) at (-3.2,-2.4) {};
  \node (n5) at (-8,-6) {};
  \node (n7) at (5,-4) [label = {below:$\partial D$}]{};

  \node (n4) at (0,-2.4) [point, label = {below:$x_k$}]{};
  \node (n5) at (3.6,-2.4) [point, label = {below:$y_k^2$}]{};
  \node (n7) at (-3.68,-2.76) [point, label = {below:$y_k^1$}]{};

  % triangle edges: connect the vertices, and leave a node at the midpoint
  \draw[-] (8,-3) -- (0,3);
  \draw[-] (-8,-3) -- (0,3);
  \draw[-] (-3.2,-2.4) -- (-8,-6);
  \draw[-] (-3.2,-2.4) -- (5,-2.4) -- (5,-4);
\end{tikzpicture}
\caption{}
\label{fig1}
\end{figure}

\begin{figure}[H]
\centering
\begin{tikzpicture}
  [
    scale=0.75,
    %>=stealth,
    point/.style = {draw, circle,  fill = black, inner sep = 1pt},
    %dot/.style   = {draw, circle,  fill = black, inner sep = .2pt},
  ]

  % the circle
  \def\rad{2.4}
  \node (o1) at (0,0) [point, label = {below:$z(x_k,\theta_k)$}]{};
  \draw (o1) circle (\rad);

  % triangle nodes: just points on the circle
  \node (n1) at (0,3) [label = {left:$\partial J_k$}]{};
  \node (n2) at (-4,0) {};
  \node (n3) at (4,0) [label = {below right:$\partial A$}] {};
  \node (n4) at (0,-2.4) [point, label = {below:$x_k$}]{};
  \node (n7) at (1.5,-2.4) [label = {below:$\partial D$}]{};

  \node (n5) at (-1.44,1.92) [point, label = {below:$x_k^1$}]{};
  \node (n6) at (1.44,1.92) [point, label = {below:$x_k^2$}]{};
  \node (n8) at (0,3) [point, label = {above:$x_k^3$}]{};
  \node (n9) at (0.3,3) [label = {right:$\zeta_k$}]{};

  % triangle edges: connect the vertices, and leave a node at the midpoint
  \draw[-] (n3) -- (0,3);
  \draw[-] (n2) -- (0,3);
  \draw[-] (-2,-2.4) -- (1.5,-2.4);
  \draw[dashed] (n5) -- (n6);
  \draw[dashed] (n8) -- (1.44,4.08);
  \draw[-] (0.24,2.82) arc (-37:37:0.3);

\end{tikzpicture}
\caption{}
\label{fig2}
\end{figure}

Now we estimate $|J_k|$, we consider the two tangent lines of $B_{\rho(x_k,\theta_k)}(z(x_k,\theta_k))$ at $x_k^1$ and $x_k^2$ and the edge $\overline{x_k^1x_k^2}$. These three lines forms a triangle $J_k'$. Since $A$ is convex and $B_{\rho(x_k,\theta_k)}(z(x_k,\theta_k))$ is tangent to $A$, we have $J_k \subset J_k'$. When $\zeta_k \leq \frac{\pi}{2}$, we can compute the area of the triangle $J_k'$
\begin{align*}
    |J_k'| = 2\rho^2(x_k,\theta_k) \sin^2\Big(\frac{\zeta_k}{2}\Big) \tan\Big(\frac{\zeta_k}{2}\Big)
           \leq \frac{1}{2} \rho^2(x_k,\theta_k) \zeta_k^3.
\end{align*}
When $\zeta_k > \frac{\pi}{2}$, we bound the area of the triangle $J_k'$ by the diameter of the set A
\begin{align*}
    |J_k'| = 2\diam A \rho(x_k,\theta_k) \sin\Big(\frac{\zeta_k}{2}\Big) \leq 2\diam A \rho(x_k,\theta_k).
\end{align*}
Note that $\sum_{k \in \Lambda} \zeta_k \leq 2\pi$ and $\zeta_k \geq 0$, then there are at most four connected components with index $k \in \Lambda$ such that $\zeta_k > \frac{\pi}{2}$. Therefore we have
\begin{align*}
    |J| &= \sum_{k \in \Lambda} |J_k| \\
        &\leq \sum_{k \in \Lambda, \zeta_k > \frac{\pi}{2}} |J_k'| + \sum_{k \in \Lambda, \zeta_k \leq \frac{\pi}{2}} |J_k'| \\
        &\leq 8 \diam A \sup_{k \in \Lambda}\rho(x_k,\theta_k) + \frac{1}{2} \sum_{k \in \Lambda, \zeta_k \leq \frac{\pi}{2}} \rho^2(x_k,\theta_k) \zeta_k^3 \\
        &\leq \Big( 8 + \frac{\pi^3}{4} \Big) R \diam A
\end{align*}
Now insert the last inequality and \Cref{basic_convex} into \eqref{geometry_lem_ann_eq4}. We have
\begin{align*}
    |A \backslash \bar{D}| \leq \Big( 8 + 3\pi + \frac{\pi^3}{4} \Big) R \diam A.
\end{align*}
This concludes the proof.
\end{proof}

% !TEX root = P5_paper1.tex

\section{H\"older continuity of the minimizer} \label{regularity}

In this section we prove the global minimizer $\omega$ is H\"older continuous.

\begin{theorem} \label{regularity_degenerate}
Given a bounded convex open domain $\Omega \subset \R^2$ with $C^{2,\alpha/4}$ boundary, $\alpha \in (0,1]$. Let $\omega_0 \in C^\alpha(\bar{\Omega})$ satisfy $\inf_{x \in \Omega} \omega_0(x) >0$ or $\sup_{x \in \Omega} \omega_0(x) <0$. Then the global minimizer $\omega$ of the functional $E$ \Cref{main_thm_intro} is in $C^{\alpha/4}(\bar{\Omega})$.
\end{theorem}

The key to show $f = (\omega^*)^{-1} \circ \psi^*$ is H\"older continuous. We prove $(\omega^*)^{-1}$ and $\psi^*$ are H\"older continuous in \Cref{diff_lemma_1} and \Cref{holder_half}. This is somehow difficult because a priori we have little information on the set of stagnation points of $\psi$. In fact, it has at least one stagnation point. $\psi^*$ may be irregular close to stagnation points. Here, we prove that $\psi^*$ is locally $C^{1/2}$ when locally $\omega \neq 0$, by a scaling argument which is robust for stagnation points.

\begin{lemma} \label{diff_lemma_1}
For bounded convex open domain $\Omega$ and any vorticity function $\omega \in C^{\alpha}(\xoverline{\Omega}), \, \alpha \in (0,1]$, there is a unique function $(\omega^*)^{-1}: [0, \mu(\Omega)] \rightarrow [\min \omega, \max \omega]$ such that $(\omega^*)^{-1}$ is continuous, non-decreasing with $(\omega^*)^{-1} \circ \omega^* = id$, where $id$ is the identity map. Moreover, $(\omega^*)^{-1} \in C^{\alpha/2}( [0, \mu(\Omega)] )$.
\end{lemma}

\begin{proof}
Fix $t \in (\min \omega, \max \omega]$. Taking the H\"older difference quotient of $(\omega^*)^{-1}$ at $\omega^*(t)$, for $s$ close to $t^-$, we consider
\begin{align} \label{diff_lemma_1_eq1}
  \frac{t-s}{|\omega^*(t)-\omega^*(s)|^{\alpha/2}}
    = \frac{t-s}{ | \{ x \in \Omega : s \leq \omega < t \} |^{\alpha/2} }.
\end{align}
Since $\|\omega\|_{C^\alpha} \leq C$ for some $C>0$, for some $x \in \Omega$ with $\omega(x)= \frac{s+t}{2}$, we have
\[
    \Big\{y \in \Omega \mid |x-y| < C^{-1} \Big(\frac{t-s}{2}\Big)^{\frac{1}{\alpha}} \Big\} \subset \{ x \in \Omega : s \leq \omega < t \}.
\]
Since $\Omega$ is a bounded convex open domain, there exists $\delta := \delta(\Omega)>0$, such that
\[
    \Big| \Big\{y \in \Omega \mid |x-y| < C^{-1} \Big(\frac{t-s}{2}\Big)^{\frac{1}{\alpha}} \Big\} \Big|
        \geq C^{-2} \delta \Big(\frac{t-s}{2}\Big)^{\frac{2}{\alpha}}
\]
Therefore, the difference quotient \eqref{diff_lemma_1_eq1} is bounded from above. The proof for $t=\min \omega$ is analogues.
\end{proof}

\begin{lemma} \label{holder_half}
Let $\Omega,\, \psi$ and $f$ be as in \Cref{levelset_conv}. Then $\psi^*: [\min \psi, 0] \rightarrow [0, |\Omega|]$ belongs to $C^{1/2}([\min \psi, 0])$.
\end{lemma}

\begin{proof}
We prove $\psi^*$ is in $C^{1/2}$ by contradiction. Assume $\psi^*$ is not in $C^{1/2}$ around $t \in [\min \psi, 0]$, then there exists $\{t_n\}_{n \geq 1}$ converging to $t$ and $\{h_n\}_{n \geq 1} \subset \R^+$ converging to $0$ such that
\begin{align} \label{holder_half_eq1}
    \frac{\psi^*(t_n+h_n)-\psi^*(t_n)} {\sqrt{h_n}} = C_n \xrightarrow{n \rightarrow \infty} \infty.
\end{align}
Since $\{\psi \leq t_n+h_n\}$ and $\{\psi \leq t_n\}$ are both convex sets, by \Cref{geometry_lem_ann}, there exists $B_{r_n}(x_n) \subset \{t_n \leq \psi \leq t_n+h_n\}$ with
\begin{align} \label{holder_half_eq3}
    r_n \geq \frac{\epsilon(\psi^*(t_n+h_n)-\psi^*(t_n))}{\diam(\Omega)}.
\end{align}
Now we zoom in these small balls. Define $\phi_n(x) := r_n^{-2} (\psi(r_n x+x_n) - \psi(x_n))$ for $x \in B_1(0)$, then $\Delta \phi_n \in [f(t_n),f(t_n+h_n)]$. Therefore, $\Delta \phi_n$ converges to the constant function $f(t)$ in $L^\infty$ as $n \rightarrow \infty$.

Now we prove $\phi_n \rightarrow 0$ in $L^\infty$. Indeed, with \eqref{holder_half_eq1} and \eqref{holder_half_eq3}, we have
\begin{align*}
    |\phi_n(x)| &\leq \Big| \frac{\diam(\Omega)^2 (\psi(r_n x+x_n) - \psi(x_n))}
                            {\epsilon^2(\psi^*(t_n+h_n)-\psi^*(t_n))^2} \Big|\\
                &\leq \frac{\diam(\Omega)^2 h_n}{\epsilon^2(\psi^*(t_n+h_n)-\psi^*(t_n))^2} \\
                &\leq \frac{\diam(\Omega)^2}{\epsilon^2C_n^2} \rightarrow 0.
\end{align*}

Decompose the sequence $\{\phi_n\}_{n \geq 1}$ into a harmonic part and a zero-trace part. Namely, let $\phi_n = \varphi_n + \xi_n$, where $\xi_n$ is a harmonic function with $\xi_n = \phi_n$ on $\partial B_1(0)$, and $\Delta \varphi_n = \Delta \phi_n$ with $\varphi_n = 0$ on $\partial B_1(0)$. Now by maximum principle, $\xi_n \rightarrow 0$ in $L^\infty$, and thus $\varphi_n \rightarrow 0$ in $L^\infty$.

On the other hand, since $\Delta \varphi_n = \Delta \phi_n$ converges to the constant function $f(t)$ in $L^\infty$ and thus $\varphi_n \rightarrow \varphi$ in $W^{1,q}$ for any $q \in (1, \infty)$. Here, $\varphi$ is the function satisfying $\Delta \varphi = f(t) >0$ and $\varphi = 0$ on $\partial B_1(0)$, which gives a contradiction.

\end{proof}

Now we can prove \Cref{regularity_degenerate}.

\begin{proof}[Proof of \Cref{regularity_degenerate}]
By \Cref{burton_thm} of Burton, we know we have a solution $\psi$ to \eqref{semilinear_Euler} with a non-decreasing function $f \in L^\infty$ and its vorticity is given by $\omega = f(\psi)$. We also have $\psi \in C^{1,\beta}$ for some $\beta>0$. Since $\omega_0 \in C^\alpha(\xoverline{\Omega})$, we have $(\omega^*)^{-1} = (\omega_0^*)^{-1} \in C^{\alpha/2}$ by \Cref{diff_lemma_1}. \par

We claim that $f$ is continuous. Indeed, if $f$ is not continuous in $(\min \psi, \max \psi)$, then there exist $s \in (\min \psi, \max \psi)$, $t_1,t_2 \in [\min \omega, \max \omega]$ with $t_1 < t_2$ such that $f(s^-)=t_1$ and $f(s^+)=t_2$. Because $f$ is non-decreasing, this contradicts with the continuity of $\omega_0$ and $\omega^* = \omega_0^*$. If $f$ is not continuous at $\min \psi$ or $\max \psi$, then $\{\psi = \min \psi\}$ or $\{\psi = \max \psi\}$ has measure zero, otherwise it contradicts with the continuity of $\omega_0$ and $\omega^* = \omega_0^*$. Hence $f(\psi)$ admits a continuous representation in $[\min \psi, \max \psi]$ and $f$ is also continuous. Now we apply \Cref{levelset_conv} and \Cref{holder_half} to deduce $\psi^* \in C^{1/2}([\min \psi, \max \psi])$.

From \eqref{semilinear_Euler}, we have $\omega^* = (f \circ \psi)^*$. Since $\omega^*: [\min \omega, \max \omega] \rightarrow [0, \mu(\Omega)]$ is strictly increasing, \Cref{diff_lemma_1} shows that it admits an non-decreasing and continuous left inverse denoted by $(\omega^*)^{-1}$, which gives $(\omega^*)^{-1} = \big[ (f \circ \psi)^* \big]^{-1}$. Also, we claim that $\psi^*$ admits an inverse denoted by $(\psi^*)^{-1}$. Indeed, note that $\omega>0$. By \Cref{levelset_conv}, $\{\psi = c\}$ has measure zero and $\{\psi \leq c\}$ is convex for any $c \in [\min \psi, \max \psi]$, so $\psi^*$ is strictly increasing and continuous. Since $f$ is continuous and non-decreasing, we have $(\omega^*)^{-1} = f \circ (\psi^*)^{-1}$, and hence $f = (\omega^*)^{-1} \circ \psi^*$ and $f \in C^{\alpha/4} ([\min \psi, \max \psi])$. Then $\omega = f(\psi) \in C^{\alpha/4}(\bar{\Omega})$ and by Schauder theory $\psi \in C^{2,\alpha/4}(\bar{\Omega})$.
\end{proof}

% !TEX root = P5_paper1.tex

\section{The structure of the critical points} \label{critical_points}

In this section, we prove the set of stagnation points is exactly $\{\psi = \min \psi\}$. Thus, by the results Nadirashvili showed in \cite{MR3067825}, the level sets $\{\psi = s\}$ for any $s<0$ are analytic curves or a single point. The smoothness of $\omega$ in any compact subset of $\Omega \backslash \{\psi = \min \psi\}$ also follows from an iterative argument in Theorem 1.5 of \cite{MR3067825}. Therefore, (1) and (4) in \Cref{main_thm_intro} follow from \Cref{regularity_degenerate} the following proposition.

\begin{proposition} \label{no_critical}
Given $\Omega,\, \omega_0$ and $\omega$ as in \Cref{main_thm_intro}, let $\psi = (\Delta)^{-1} \omega$ be the global minimizer of the kinetic energy functional $E$ in \eqref{kinetic_energy}. Then $|\nabla \psi(x)| \neq 0$ for any $x \in \bar{\Omega}$ with $\psi(x) \neq \min \psi$. Therefore, the set of stagnation points is exactly $\{\psi = \min \psi\}$ which is a single point or a line segment.
\end{proposition}

\begin{proof}
Assume there exists $x_0 \in \Omega$ with $|\nabla \psi(x_0)| = 0$ and $\psi(x_0) = c \in (\min \psi, \max \psi)$. Let $\Delta \psi(x_0) = a > 0$. Since $\omega \in C^{\alpha/4}$, locally after a conformal change of coordinate, $\psi(x^1,x^2) = a_1 (x^1)^2 + a_2 (x^2)^2 + o(x^1,x^2) + c$ with $a=a_1+a_2$ and $|o(x^1,x^2)| \leq C(|x^1|^{2+\alpha/4} + |x^2|^{2+\alpha/4})$. About the local behavior of $\psi$, we have the following situations:
\begin{enumerate}[leftmargin=*,label=\textup{(\arabic*)},align=left]
    \item $a_1a_2>0$ is impossible. Indeed, in this case $\psi$ achieves a local strict minimum at $x=x_0$.
    \item $a_1a_2<0$ is also impossible. Indeed, in this case the level set is locally homeomorphic to the set $\{(x^1,x^2) \in \R^2 \mid |x^1| = |x^2| \leq 1\}$. However, the level set should be a convex curve by \Cref{levelset_conv}.
    \item $a_1a_2=0$ with $a_1>0$, then $\{\psi < c\} \cap \{x^2=0\} = \varnothing$. Indeed, if $\{\psi < c\} \cap \{x^2=0\} \neq \varnothing$, for any $x \in \{\psi < c\} \cap \{x^2=0\}$, $\psi(x) = a_1 (x^1)^2 + o(x^1,0) + c > c$. Moreover, for any $x \in \{\psi=c\}$ around $x_0$, we have $a_1 (x^1)^2 + o(x^1,x^2) = 0$, and hence
    \begin{align*}
    	a_1 (x^1)^2 \leq C(|x^1|^{2+\alpha/4} + |x^2|^{2+\alpha/4}).
    \end{align*}
    Then there exists $\delta>0$ small enough, such that for $|x^1| \leq \delta$, we have
    \begin{align} \label{no_critical_eq5}
    	\sqrt{\frac{a_1}{2C}} x^1 \leq |x^2|^{1+\alpha/8}.
    \end{align}
    Now we know that \eqref{no_critical_eq5} encloses a cusp-like region. The only possibility for the convex curve in this region containg the origin is the line $\{x^1=0\}$. This means the level set $\{\psi=c\}$ is a single line and it contradicts with \Cref{levelset_conv}. This situation is also impossible.
\end{enumerate}
Now we have exclude the possibility of $|\nabla \psi(x_0)| = 0$ for $\psi(x_0) = c \in (\min \psi, \max \psi)$. The case $c = \max \psi$ is only achieved on the boundary. Because the boundary is $C^{2,\alpha/4}$, we can apply Hopf boundary lemma to deduce $|\nabla \psi(x)| \neq 0$. By \Cref{levelset_conv}, we know $\{\psi = \min \psi\}$ is a single point or a line segment.
\end{proof}

% !TEX root = P5_paper1.tex

\section{Constraints on level set topology} \label{counterexamples}

In this section, we show there exists a unique global minimizer in $L^\infty$-closure of $\mathcal{O}_{\omega_0}$ in \Cref{counterexample_2}, when $\omega_0$ is H\"older continuous and has nice level set topology. We also prove these conditions are sharp in \Cref{counterexample_1} and \Cref{counterexample_3}.

\begin{theorem} \label{counterexample_2}
Given a bounded convex domain $\Omega \subset \R^2$, $\alpha \in (0,1]$. Let $\omega_0 \in C^\alpha(\bar{\Omega})$ satisfy $\omega_0 \geq 0$. Suppose $\omega_0$ is constant on the boundary and the set $\{\omega_0 < s\}$ is simply connected for any $s \in (\inf \omega_0, \sup \omega_0)$, then $E$ has a global minimizer given by $\omega$ in \Cref{burton_thm} in $L^\infty$-closure of $\mathcal{O}_{\omega_0}$.
\end{theorem}

\begin{theorem} \label{counterexample_1}
Let $\Omega,\,\alpha$ as in \Cref{counterexample_2}. Let $\omega_0 \in C^\alpha(\bar{\Omega})$ satisfy $\omega_0 \geq 0$. Assume there exists some $s \in (\inf \omega_0, \sup \omega_0)$ such that the set $\{\omega_0 < s\}$ is not simply connected or $\omega_0$ is not constant on the boundary, then $E$ has no global minimizer in $L^\infty$-closure of $\mathcal{O}_{\omega_0}$.
\end{theorem}

\begin{theorem} \label{counterexample_3}
Let $\Omega,\,\alpha$ as in \Cref{counterexample_2}. There exists $\omega_0 \in L^\infty(\Omega)$ with $\omega_0 \geq 0$ satisfying the following conditions. The set $\{\omega_0 < s\}$ is simply connected for any $s \in (\inf \omega_0, \sup \omega_0)$, and $\omega_0$ is constant on the boundary. Moreover, $E$ has no global minimizer in $L^\infty$-closure of $\mathcal{O}_{\omega_0}$.
\end{theorem}

These results need \Cref{action2}. \Cref{action2} is a lemma about the deformation of volume-preserving diffeomorphisms. For any simply connected open sets $G, D \subset \R^2$ with $\bar{G}, \bar{D} \subset B_1(0)$, there exists a volume-preserving homeomorphism $\eta$ with $\eta (\partial B_1(0)) = \partial B_1(0)$ and $\eta(G)=D$, if and only if $D$ and $G$ have the same measure. \Cref{action2} is a corollary of this fact where we need local diffeomorphisms instead of global homeomorphisms. We omit the proof of \Cref{action2}.

\begin{lemma} \label{action2}
Given a bounded connected open set $A$, for any open set $G$ with $\bar{G} \subset A$ which is homeomorphic to an annulus, define open sets $G_i$ and $G_o$ by the unique decomposition $A = G_i \cup \bar{G} \cup G_o$ such that:
\begin{enumerate}[leftmargin=*,label=\textup{(\arabic*)},align=left]
    \item The intersection of any two of $G_i, G$ and $G_o$ is empty;
    \item $G_i$ is simply connected. $G_o$ is homeomorphic to an annulus.
\end{enumerate}
Given another open set $D$ with $\bar{D} \subset A$ which is homeomorphic to an annulus, define $D_i$ and $D_o$ similarly. Suppose we have $|G_i| < |D_i|$ and $|G_o| < |D_o|$. Then there exists a volume-preserving diffeomorphism $\eta:A \rightarrow A$ such that $\eta(D) \subset G$ and $\eta = id$ in a neighborhood of $\partial A$.
\end{lemma}

The following remark is useful in the proof of \Cref{ce_l1}.

\begin{remark} \label{rmk_action2}
If we only know $D$ and $G$ are open subsets of $A$ and $|D| \leq |G|$, without any topological assumptions, we can still find an open set $D' \subset D$ with $|D|-|D'| = \epsilon$ for any $\epsilon \in (0, |D|)$, such that there exists a volume-preserving diffeomorphism $\eta:A \rightarrow A$ with $\eta(D') \subset G$ and $\eta = id$ in a neighborhood of $\partial A$. This is done by cutting $G$ into $G' \subset G$ such that $|G| - |G'| = \frac{\epsilon}{2}$, every connected component of $G'$ is simply connected, $G'$ has finitely many connected components and $\partial G'$ has finite $1$-dimensional Hausdorff measure. We also cut $D$ into $D'$ similarly with an additional condition: they fit into the connected components of $G'$. Then we apply the simple observation above for each connected component of $D'$.
\end{remark}

Now we prove \Cref{counterexample_2}, \Cref{counterexample_1} and \Cref{counterexample_3}.

\begin{proof}[Proof of \Cref{counterexample_2}]
It suffices to show for any $\delta>0$, we can find $\omega' \in \mathcal{O}_{\omega_0}$ such that $\|\omega'-\omega\|_{L^\infty} \leq \delta$. This is done in $N:=\lceil \frac{2(\sup \omega_0 - \inf \omega_0)}{\delta} \rceil$ steps. Note that the level set topology does not change after we apply a volume-preserving diffeomorphism. Namely, for any volume-preserving diffeomorphism $\eta$, $\omega_0 \circ \eta$ is constant on $\partial \Omega$ and $\{\omega_0 \circ \eta < s\}$ is simply connected for any $s \in (\inf \omega_0 \circ \eta, \sup \omega_0 \circ \eta)$.
\par

In the first step, we apply \Cref{action2} with $A= \Omega$, $G=\{\sup \omega - \frac{2\delta}{3} < \omega < \sup \omega - \frac{\delta}{3}\}$, $D=\{\sup \omega_0 - \frac{5\delta}{8} < \omega_0 < \sup \omega_0 - \frac{3\delta}{8}\}$, then we get $\eta_1$. By \Cref{burton_thm} and \Cref{levelset_conv}, all asumptions in \Cref{action2} are satisfied. In the $k$-th step, we apply \Cref{action2} with
\begin{align*}
	A&= \Big\{\omega_0 \circ \eta_{k-1}\circ \cdots \circ \eta_1 < \sup \omega_0 -\frac{(k-1)\delta}{2} -\frac{\delta}{8} \Big\},\\
	G&= \Big\{\sup \omega - \frac{k\delta}{2}-\frac{\delta}{6} < \omega < \sup \omega - \frac{k\delta}{2} + \frac{\delta}{6} \Big\},\\
	D&= \Big\{\sup \omega_0 - \frac{k\delta}{2}-\frac{\delta}{8} < \omega_0 < \sup \omega_0 - \frac{k\delta}{2} + \frac{\delta}{8} \Big\},\\
\end{align*}
then we get $\eta_k$. Here, $G$ and $D$ are always homeomorphic to a ring, and $A$ is simply connected. After $N$ steps, we have
\[
	\|\omega - \omega_0 \circ \eta_{N} \circ \eta_{N-1} \circ \cdots \circ \eta_1 \|_{L^\infty} \leq \delta.
\]
Letting $\omega' := \omega_0 \circ \eta_{N} \circ \eta_{N-1} \circ \cdots \circ \eta_1 \|_{L^\infty} \leq \delta$ concludes our proof.
\end{proof}

\begin{proof}[Proof of \Cref{counterexample_1}]
Assume $E$ has a global minimizer $\omega'$ in $L^\infty$-closure of $\mathcal{O}_{\omega_0}$. This minimizer is necessarily $\omega$ in \Cref{burton_thm}. Indeed, by \Cref{burton_thm}, there exist $\{\omega_k\}_{k \geq 1} \subset \mathcal{O}_{\omega_0}$ with $\omega_k \rightarrow \omega$ in $L^q$ for any $q \in [1,\infty)$, then
\[
	\int_\Omega |\nabla \Delta^{-1} \omega'|^2 dx \geq \int_\Omega |\nabla \Delta^{-1} \omega|^2 dx
	= \lim_{k \rightarrow \infty} \int_\Omega |\nabla \Delta^{-1} \omega_k|^2 dx \geq \int_\Omega |\nabla \Delta^{-1} \omega'|^2 dx.
\]
Since $\omega$ is the unique minimizer of $E$ in the $L^q$-strong closure of $\mathcal{O}_{\omega_0}$ for any $q \in [1,\infty)$, then $\omega' = \omega$. \par

Now we show $\omega$ is not in the the $L^\infty$-strong closure of $\mathcal{O}_{\omega_0}$. If the set $\{\omega_0 < s\}$ is connected but not simply connected for some $s \in (\inf \omega_0, \sup \omega_0)$, we have two possibilities:
\begin{enumerate}[leftmargin=*,label=\textup{(\arabic*)},align=left]
    \item There exists $\epsilon$ such that the set $\{s-\epsilon < \omega_0 < s\}$ contains at least two connected components $D_1, D_2$ with $\dist(D_1,D_2)>0$. However, the level set $\{s-\epsilon < \omega < s\}$ is a convex ring, hence a connected set. By some elementary topological argument, we can show for any $v$ in the the $L^\infty$-strong closure of $\mathcal{O}_{\omega_0}$, $\|v-\omega\|_{L^\infty}> \frac{\epsilon}{2}$. Therefore, $\omega$ is not in the the $L^\infty$-strong closure of $\mathcal{O}_{\omega_0}$.
    \item $\partial \Omega \subset \partial \{\omega_0 < s\}$. Then $\{\omega_0>s\}$ is in the interior of $\Omega$. From \Cref{levelset_conv} and the fact that $f$ is nondecreasing, there exists a neighborhood $O$ of the boundary $\partial \Omega$, such that $\omega > \frac{\sup \omega + s}{2}$ in $O$. Then for any $v$ in the the $L^\infty$-strong closure of $\mathcal{O}_{\omega_0}$, $\|v-\omega\|_{L^\infty} \geq \frac{\sup \omega - s}{2}$.
\end{enumerate}
 The case where the set $\{\omega_0 < s\}$ is not connected follows from the same argument as (1). If $\omega_0$ is not constant on the boundary, the same conclusion follows from an argument similar to (2).
\end{proof}

\begin{proof}[Proof of \Cref{counterexample_3}]
As in the proof of \Cref{counterexample_1}, assume $E$ has a global minimizer in $L^\infty$-closure of $\mathcal{O}_{\omega_0}$, then this global minimizer must be $\omega$ in \Cref{burton_thm}. We construct the following vortex patch supported on $\Omega' \subset \Omega$ with $\bar{\Omega}' \subset \Omega$ and $\partial \Omega'$ is not locally Lipschitz (for example, $\partial \Omega'$ contains a cusp-like corner):
\begin{align}
	\omega_0(x)=
	\begin{cases}
		1, \quad x \in \Omega', \\
		0, \quad x \in \Omega \backslash \Omega'.
	\end{cases}
\end{align}
Note that $\omega$ has convex level sets by \Cref{burton_thm} and \Cref{levelset_conv}, and any volume-preserving diffeomorphism maps convex level curves to locally Lipschitz curves. Then for any $\omega' \in \mathcal{O}_{\omega_0}$, we have $\|\omega'-\omega\|_{L^\infty}=1$. This gives a contradiction.
\end{proof}

\appendix
% !TEX root = P5_paper1.tex

\section{A counterexample} \label{counterexample_naive_sec}

In this section, we give a counterexample to show that generally there may not exist smooth minimisers for the kinetic energy functional $E$ in an orbit $\mathcal{O}_\omega$, even for smooth $\omega_0$ with nice level set topology.

\begin{theorem} \label{counterexample_naive}
Let $\Omega = B_1$ be the unit disk in $\R^2$, then there exists a nonnegative smooth function $\omega_0$ with a single critical point and $\omega_0 |_{\partial B_1} = c$ for some $c>0$, such that $E$ admits no smooth minimiser in the orbit $\mathcal{O}_{\omega_0}$.
\end{theorem}

To prove \Cref{counterexample_naive}, we need and a variant of P\'olya–Szeg\H{o} inequality for fourth order differential operators. This have been proved in \cite{MR1135975}. It can also be deduced from classical symmetrization techniques, see Chapter 2 and Lemma 11.3.2 in \cite{MR2251558}. We also need \Cref{ce_l1} for which a more general version might be known elsewhere. Since we cannot locate a direct and precise reference, we still give a concise proof.

\begin{lemma} \label{PS_inequality}
For any nonnegative function $u \in W^{1,q}(B^1)$ with $B^1 \subset \R^n$, $n \in \N^+$ and $1 \leq q \leq \infty$, define its symmetric increasing rearrangement
\[ \tilde{u}(x) = \sup_{\mu (\{u \leq s\}) \leq \mu(B_{|x|})} s, \]
then
\[ \int_{\R^n} |\nabla \Delta^{-1} \tilde{u}|^q dx \leq \int_{\R^n} |\nabla \Delta^{-1} u|^q dx. \]
\end{lemma}

\begin{lemma} \label{ce_l1}
Given a bounded domain $\Omega \subset \R^2$, let $\omega_0,\, \omega \in C(\bar{\Omega})$ with $\omega_0^* = \omega^*$. Then there exist $\{\omega_k\}_{k \geq 1} \subset \mathcal{O}_{\omega_0}$ with $\omega_k \rightarrow \omega$ in $L^q$ as $k \rightarrow \infty$ for any $q \in [1,\infty)$.
\end{lemma}

\begin{proof}[Proof of \Cref{counterexample_naive}]
We design $\omega_0$ such that $\omega_0(x,y) = 1+2(x^2+y^4)$ in $B_{3/4}$. Moreover, $\omega_0$ should satisfy all properties mentioned in \Cref{counterexample_naive}. By elementary computation, we have for any $x^2+y^2 < 9/16$,
\[ \tilde{\omega}_0(x,y) = 1 + 2\Big(\frac{\pi}{\mu(\{(x_1,y_1) \,|\, x_1^2+y_1^4 \leq 1\})}\Big)^{4/3} (x^2+y^2)^{4/3}.  \]
We can use \Cref{PS_inequality} to deduce
\[ \int_{B_1} |\nabla \tilde{\omega}_0|^2 dx \leq \int_{B_1} |\nabla \omega|^2 dx 
	\quad \text{for any } \omega \in \mathcal{O}_{\omega_0}. \]
\par

To conclude the proof, it suffices to show there are $\{\omega_k\}_{k \in \N} \subset \mathcal{O}_{\omega_0}$ such that
\begin{align} \label{counterexample_naive_eq3} 
	\lim_{k \rightarrow \infty} \int_{B_1} |\nabla \omega_k|^2 dx
		= \int_{B_1} |\nabla \tilde{\omega}_0|^2 dx.
\end{align}
This is a consequence of \Cref{ce_l1}.
\end{proof}

\begin{proof}[Proof of \Cref{ce_l1}]
It suffices to prove that for any $\delta>0$, we can find $\omega' \in \mathcal{O}_{\omega_0}$ such that $\|\omega'-\omega\|_{L^1} \leq \delta$. \par

Define a sequence $\{a_k\}_{k \geq 1} \subset \N^+$ with $a_1 = 1$ and $a_{k+1} = 1+\sum_{l=1}^{k} a_l$ for any $k \geq 1$. Let $\kappa = \frac{\delta}{2|\Omega|}$, $N = \lceil \frac{\sup \omega_0 - \inf \omega_0}{\kappa} \rceil$ and $\epsilon = \frac{\delta}{2 \sum_{k=1}^N a_k (\sup \omega_0 - \inf \omega_0)}$. Then we construct $\omega'$ in $N$ interative steps. In the first step, consider the set $\{\omega \leq \kappa + \inf \omega\}$. If $\{\omega \leq \kappa + \inf \omega\}$ is empty, we skip this step. Otherwise, we apply \Cref{rmk_action2} with $A = \Omega$, $\epsilon$ defined above, $G = \inte \{\omega \leq \kappa + \inf \omega\}$ and $D'$ to be any subset of $\inte \{\omega_0 \leq \kappa + \inf \omega_0\}$ with $|G|-|D'| \leq \epsilon$ and satisfying the requirement in \Cref{action2}, then we have a volume-preserving diffeomorphism denoted by $\eta_1$. Then
$$ \int_{\{\omega \leq \kappa + \inf \omega\}} |\omega - \omega_0 \circ \eta_1| dx
	\leq \kappa |\{\omega \leq \kappa + \inf \omega\}| + \epsilon (\sup \omega_0 - \inf \omega_0) a_1. $$
In $k$-th step, we apply \Cref{action2} and \Cref{rmk_action2} with $A = \Omega \backslash \{\omega \leq (k-1)\kappa + \inf \omega\}$, $G = \{(k-1)\kappa + \inf \omega < \omega \leq k \kappa + \inf \omega\}$ and $D'$ to be any subset of $\{(k-1)\kappa + \inf \omega_0 < \omega_0 \leq k \kappa + \inf \omega_0\} \cap A$ with $|G|-|D'| \leq \epsilon a_k$ and satisfying the requirement in \Cref{action2}, (we skip this step if $G$ is empty) then we get $\eta_k$ with
\begin{align*}
	\int_{\{(k-1)\kappa + \inf \omega < \omega \leq k \kappa + \inf \omega\}}& |\omega - \omega_0 \circ \eta_k \circ \cdots \circ \eta_1| dx \\
	&\leq \kappa |\{(k-1)\kappa + \inf \omega < \omega < k \kappa + \inf \omega\}| + \epsilon (\sup \omega_0 - \inf \omega_0) a_k.
\end{align*}
Finally, let $\omega' = \omega_0 \circ \eta_{N} \circ \eta_{N-1} \circ \cdots \circ \eta_1$, and we have
\begin{align*}
	\int_{\Omega} |\omega-\omega'| dx \leq \kappa |\Omega| + \epsilon (\sup \omega_0 - \inf \omega_0) \sum_{k=1}^{N} a_k \leq \delta.
\end{align*}
Since $\mathcal{O}_{\omega_0}$ is uniformly bounded in $L^\infty$, we have the convergence in $L^q$ for any $q \in (1,\infty)$.
\end{proof}

%uncomment next line to change bibliography name to references
\renewcommand{\bibname}{References}
\bibliography{bibliography}
\bibliographystyle{abbrv}

\end{document}